\documentclass[12pt, leqno]{article}
\usepackage{amsmath}
\usepackage{amssymb}
\usepackage{amsthm}
\def\underset#1#2{{\mathrel{\mathop {{}_{} {#2}}\limits_{{#1}_{}}}}}
\def\upplim_#1{\underset{#1}{\overline\lim}\;}
\def\lowlim_#1{\underset{#1}{\underline\lim}\;}
\newcommand{\C}{{\mathbf{C}}}
\newcommand{\codim}{{\mathrm{codim}}\,}
\newcommand{\del}{{\partial}}
\newcommand{\delbar}{\bar{\partial}}
\newcommand{\fa}{\forall}

\newcommand{\grad}{{\mathrm{grad}}}

\newcommand{\lam}{{\lambda}}

\newcommand{\pnc}{{\mathbf{P}^n(\mathbf{C})}}
\newcommand{\ponec}{{\mathbf{P}^1(\mathbf{C})}}
\newcommand{\R}{{\mathbf{R}}}

\newcommand{\Z}{{\mathbf{Z}}}

\newcommand{\tensor}{\otimes}

\newtheorem{cor}[equation]{Corollary}

\newtheorem{ex}[equation]{\bf Example}

\newtheorem{thm}[equation]{\bf Theorem}
\newtheorem{mthm}[equation]{\bf Main Theorem}
\newtheorem{prob}[equation]{\bf Problem}

\newtheorem{rmk}[equation]{\indent {\it Remark}\rm}
\numberwithin{equation}{section}

\title{Order of Meromorphic Maps and Rationality\\ of the Image Space%
\thanks{Research supported in part by Grant-in-Aid
for Scientific Research (A) 60218790 and SFB/TR 12 (DFG).
}}
\date{\hfil}
\author{Junjiro Noguchi and J\"org Winkelmann}
\begin{document}
\parindent12pt
\baselineskip18pt
\maketitle
\begin{abstract}
Let $\iota : \C^2 \hookrightarrow S$ be a compactification
of the two dimensional complex space $\C^2$.
By making use of  Nevanlinna theoretic methods and
the classification of compact complex surfaces
K. Kodaira proved in 1971 (\cite{ko71})
that $S$ is a rational surface.
Here we deal with a more general meromorphic map $f: \C^n \to X$ into a
compact complex manifold $X$ of dimension $n$, whose differential
$df$ has generically rank $n$.
Let $\rho_f$ denote the order of $f$. We will prove that
if $\rho_f<2$, then every global symmetric
holomorphic tensor must vanish; in particular,
{\it if $\dim X=2$ and
$X$ is k\"ahler, then $X$ is a rational surface.
Without the k\"ahler condition there is no such conclusion, 
as we will show by a counter-example using a Hopf surface.}
This may be the first instance that the k\"ahler or non-k\"ahler
condition makes a difference in the value distribution theory.
\medskip

{\small
}
\end{abstract}

\section{Introduction and main results.}

Let $X$ be a compact hermitian manifold with metric form $\omega$.
Let $f:\C^n \to X$ be a meromorphic map (cf.\ \cite{no90} for
this section in general).
If the differential $df$ is generically of maximal rank,
$f$ is said to be {\it differentiably non-degenerate}.
We set
\begin{equation}
\label{alpha}
\alpha=dd^c \|z\|^2
\end{equation}
for $z=(z_j) \in \C^n$,
where $d^c=\frac{i}{4\pi}(\delbar - \del)$ and
$\|z\|^2=\sum_{j=1}^n |z_j|^2$.
We use the notation:
$$
B(r)=\{ z\in \C^n: \| z \| <r \}, \qquad S(r)=\{ z\in \C^n: \|z \|=r \} \quad (r>0).
$$
We define the {\it order function} of $f$ with respect to
$\omega$ by
\begin{equation}
\label{order}
T_f(r; \omega)=\int_1^r \frac{dt}{t^{2n-1}} \int_{B(t)}
f^*\omega \wedge \alpha^{n-1}.
\end{equation}
Then the (upper) order is defined by
$$
\rho_f=\upplim_{r \to \infty} \frac{\log T_f(r; \omega)}{\log r}.
$$
It is easy to see that $\rho_f$ is independent of the choice of
the metric (form $\omega$) on $X$.
\begin{ex}
\rm
(i) If $X=\pnc$ and $f$ is rational, then $\rho_f=0$.

(ii) Let $X$ be a compact torus. If $f: \C^n \to X$ is non-constant,
then $\rho_f \geq 2$. If $\lambda : \C^n \to X$ ($\dim X=n$)
is the universal covering map, then $\rho_{\lambda}=2$.
\end{ex}

A compact complex manifold which is bimeromorphic to $\pnc$ is called a
{\it rational variety}.
A two-dimensional compact complex manifold is called
a complex {\it surface}. If it admits a  k\"ahler metric,
it is called a k\"ahler surface.

The main result of this paper is the following:
\begin{mthm}
\label{main}
Let $X$ be a k\"ahler surface. Assume that there is a
differentiably non-degenerate meromorphic map $f: \C^2 \to X$.
If $\rho_f < 2$, then $X$ is rational.
\end{mthm}

The k\"ahler condition is necessary by the following:
\begin{thm}
\label{count}
There is a Hopf surface $S$ for which there is a differentiably
non-degenerate holomorphic map $f:\C^2 \to S$ with
$\rho_f=1$.
\end{thm}

Let $\Omega^k_X$ denote the sheaf of holomorphic $k$-forms over
a complex manifold $X$.
We denote by $S^l\Omega^k_X$ its $l$-th symmetric tensor power.
In particular, $K_X=\Omega^n_X$ \ ($n=\dim X$) denotes the
canonical bundle over $X$.

The key tool for the proof of the Main Theorem \ref{main} is:
\begin{thm}
\label{key}
Let $X$ be an $n$-dimensional compact complex manifold.
Assume that there exists a differentiably non-degenerate meromorphic map
$f: \C^m \to X\ \ (m \geq n)$ with $\rho_f <2$.
Then for arbitrary $l_k \geq 0$ with $\sum_{k=1}^n l_k >0$ 
$$
H^0(X, S^{l_1}\Omega^1_X \tensor \cdots \tensor S^{l_n}\Omega^n_X)=\{0\}.
$$
\end{thm}

\begin{rmk}
So far by our knowledge, the above theorems
are the first instance that the k\"ahler or non-k\"ahler
condition makes a difference in the value distribution theory.
\end{rmk}

\section{Proof of the Main Theorem.}
{\bf (1) Proof of Theorem \ref{key}.}  
Assume the existence of an element
$$
\tau \in
H^0(X, S^{l_1}\Omega^1_X \tensor \cdots \tensor S^{l_n}\Omega^n_X)
\setminus\{0\}.
$$
We take a hermitian metric $h$ on $X$ with
the associated form $\omega$.
There are induced hermitian metrics on the symmetric powers of the bundles $\Omega^k$ and their
tensor products which by abuse of notation are again by denoted by $h$.
Let $\|\tau\|_h$ denote the norm of $\tau$ with respect to $h$.
Then there is a constant $c_1>0$ such that
\begin{equation}
\label{3.3.3}
\|\tau\|_h \leq c_1.
\end{equation}
We denotes by $\xi_\lam$ the coefficient functions of
$f^*\tau$ with respect to the standard coordinate system
$(z_1, \ldots, z_m)$ on $\C^m$.
Since $f$ is meromorphic, $f^*\tau$ is obviously holomorphic outside the indeterminacy set $I_f$.
Because $\codim(I_f)\ge 2$ and because $f^*\tau$ is a section in a globally defined vector bundle,
it extends holomorphically to $I_f$. Thus we may regard $f^*\tau$ as being holomorphic on $\C^n$
and the $\xi_\lam$ are holomorphic as well.

We set
\begin{equation}
\label{3.3.4}
\|f^*\tau\|^2_{\C^m}=\sum_{\lam=1}^m |\xi_\lam|^2 \not\equiv 0 .
\end{equation}
We define a function $\zeta$ on $\C^m$ by
$$
f^*\omega\wedge \alpha^{m-1}=\zeta \alpha^m.
$$
Since $f$ is differentiably non-degenerate, $f^*\tau \not\equiv 0$.
By (\ref{3.3.3}) there are positive constants $c_2$ and $c_3$ such that
\begin{equation}
\label{3.3.5}
\zeta \geq c_2 \|f^*\tau\|^{2c_3}_{\C^m}.
\end{equation}
By (\ref{3.3.4}) $\|f^*\tau\|^{2c_3}_{\C^m}$ is plurisubharmonic.
Since $f^*\tau\not\equiv 0$ is holomorphic, it follows that 
$$
\int_{S(1)}\|f^*\tau\|^{2c_3}_{\C^m}\gamma=c_4>0,
$$
where
\begin{equation}
\label{gamma}
\gamma= \frac{1}{r^{2m-1}}\, d^c\|z\|^2\wedge\alpha^{m-1},
\end{equation}
induced on $S(r)$ with $r=1$.
Since
$$
\int_{S(r)}\|f^*\tau\|^{2c_3}_{\C^m} \gamma
$$
is monotone increasing in $r>0$,
we see that
$$
\int_{S(r)}\|f^*\tau\|^{2c_3}_{\C^m} d^c\|z\|^2\wedge
\alpha^{m-1}\geq c_4 r^{2m-1}, \qquad r >1.
$$
Therefore
$$
\int_{B(r)}\|f^*\tau\|^{2c_3}_{\C^m}
\alpha^{m}\geq \frac{c_4}{2m} (r^{2m}-1), \qquad r >1.
$$
We deduce from this that
\begin{align*}
T_f(r, \omega)&=
\int_1^r \frac{dt}{t^{2m-1}}\int_{B(t)} \zeta \alpha^m 
\geq c_2
\int_1^r \frac{dt}{t^{2m-1}}\int_{B(t)}
\|f^*\tau\|^{2c_3}_{\C^m} \alpha^m\\
&\geq \frac{c_2c_4}{2m}
\int_1^r \left(t-\frac{1}{t^{2m-1}}\right)dt
 = \frac{c_2c_4}{4m}r^2+C_m(r),
\end{align*}
where $C_1(r)=O(\log r)$ and $C_m(r)=O(1)$ for $m \geq 2$.
Thus,
$$
\rho_f=\upplim_{r\to\infty} \frac{\log T_f(r, \omega)}{\log r}
\geq 2.
$$
This is a contradiction. \hfill {\it Q.E.D.}

\begin{cor}
\label{1-dim}
If $X$ in Theorem \ref{key} is $1$-dimensional,
then $X$ is biholomorphic to $\ponec$.
\end{cor}

{\bf (2) Proof of the Main Theorem \ref{main}.}  
There is a fine classification theory of complex surfaces
 (cf.\ Kodaira \cite{ko71}, Barth-Peters-Van de Ven \cite{bpv}).
According to it we know the following fact, where
$b_1(X) =\dim H_1(X, \R)$ denotes the first Betti number of $X$.

\begin{thm}{\rm (Kodaira [68] Theorem 54)}
\label{class}
If a complex surface $X$ satisfies $b_1(X)=0$ and
$H^0(X, K_X^l)=\{0\}$ for all $l>0$,  then $X$ is rational.
\end{thm}

This enables us to prove Theorem \ref{main} as follows.
By Theorem \ref{key} $\dim H^0(X, \Omega^1_X)=0$.
Due to the k\"ahler assumption we have
$b_1(X)=2\dim H^0(X, \Omega^1_X)=0$.
Moreover, $H^0(X, K_X^l)=\{0\}$ for all $l>0$ again by Theorem \ref{key}.
It follows from Theorem \ref{class}
that $X$ is rational. \hfill {\it Q.E.D.}

\section{Proof of Theorem \ref{count}.}

Let $\lambda\in\C$ with $|\lambda|>1$. Then a Hopf surface $S$
is defined as the quotient of $\C^2\setminus\{(0,0)\}$
under the $\Z$-action given by $n:(x,y)\mapsto(\lambda^n x,\lambda^n y)$.
Such a surface $S$ is known to be diffeomorphic to $S^1\times S^3$.
As a consequence $b_1(S)=1$ and $S$ is {\it not k\"ahler}.

Now
\[
\omega= \frac{i}{2\pi} \cdot
\frac{dx\wedge d\bar{x} + d y \wedge d\bar{y}}{|x|^2+|y|^2}
= \frac{dd^c||(x,y)||^2}{||(x,y)||^2}
\]
is a positive $(1,1)$-form on $\C^2\setminus\{(0,0)\}$
which is invariant under the above given $\Z$-action.
Therefore it induces a positive $(1,1)$-form on the quotient
surface $S$ which by abuse of notation is again denoted by
$\omega$.

Let $\alpha$ and $\gamma$ be as in \eqref{alpha} and \eqref{gamma},
respectively.
We claim that the holomorphic map $f:\C^2\to S$ induced by
\[
(z,w)\mapsto (z,1+zw)
\]
is of order $1$.
By definition this means
$$
\rho_f=\upplim_{r \to \infty}
\frac{\log T_f(r,\omega)}{\log r} = 1,
$$
i.e.,
\[
\upplim_{r \to \infty} \frac{1}{\log r}\log 
\int_1^r\frac{dt}{t^3}\int_{B(t)}f^*\omega \wedge\alpha = 1 .
\]
Note that
\[
f^*\omega \wedge\alpha=\frac{1+|z|^2+|w|^2}{2(|z|^2+|1+zw|^2)}\alpha^2.
\]
We define
\[
I_r=\int_{S(r)}\frac{r^2}{|z|^2+|1+zw|^2}dV,\qquad
r=\|(z,w)\|.
\]
Here $dV$ is the euclidean volume element on $S(r)$,
and therefore a constant multiple of $r^3\gamma$.
It is sufficient to show 
\begin{equation}
\label{claim}
I_r=O(r^{2+\varepsilon}), \quad \fa\varepsilon>0,
\quad \mbox{ and }\quad r^2=O(I_r).
\end{equation}
Indeed, assume that this holds.
Because of
$\lim_{r\to\infty}\frac{1+r^2}{r^2}=1$, 
\eqref{claim} is equivalent to the assertion
\[
I'_r=O(r^{1+\varepsilon}),
\quad \mbox{ and }\quad r^2=O(I'_r).
\]
with
\[
I'_r= \int_{S(r)}\frac{1+r^2}{|z|^2+|1+zw|^2}dV.
\]
From this we first obtain
\[
\int_{B(r)}\frac{1+r^2}{|z|^2+|1+zw|^2}\alpha^2=O\left(\int^r I'_r dr
\right)=O(r^{3+\varepsilon}), \quad \fa \varepsilon>0,
\]
implying
\[
T_f(r)=\frac{1}{2}\int_1 ^r \frac{dt}{t^3}\int_{B(r)}
\frac{1+r^2}{|z|^2+|1+zw|^2}\alpha^2=O(r^{1+\varepsilon}),
\quad \forall \varepsilon>0,
\]
and
\[
\rho_f=\upplim_{r \to infty} \frac{\log T_f(r)}{\log r}\leq 1 .
\]
In the same way from the second estimate of \eqref{claim}
we get the opposite estimate
$\rho_f \geq 1$, and therefore $\rho_f=1$.
Hence it suffices to show \eqref{claim}.

We define 
\[
\eta=\frac{r^2}{|z|^2+|1+zw|^2}.
\]
Thus we have to show
\[
I_r=\int_{S(r)}\eta dV=O(r^{2+\varepsilon}).
\]
We set 
$$
\eta=\frac{r^2}{\phi(z,w)}, \qquad
\phi(z, w)=|z|^2+|1+zw|^2.
$$
\subsection{Geometric estimates.}

For $(z, w)\in S(r)$ let $\theta \in[0,2\pi)$ such that $e^{i\theta}|zw|=zw$.
Let $K>0$, $-\infty < \lambda<1$ and $\mu\geq 0$.
We set
$$
\Omega_{K,\lambda,\mu}=\{(z, w)\in S(r):
|z|\leq Kr^\lambda, \: |\sin \theta| \leq r^{-\mu}\}.
$$
We need some volume estimates.

First we note that $(\sin \theta)/\theta\ge 2/\pi$ for all
$\theta\in[0,\pi/2]$,
because $\sin$ is concave on $[0,\pi/2]$.
It follows that for every $C\in]0,1]$ we have the following bound
for the Lebesgue measure:
\begin{equation}
\label{sinvol}
\mathrm{vol}\left(\{\theta\in[0,2\pi]:
|\sin \theta|\le C\}\right)\ \leq\ 4(C\pi/2)=2C\pi.
\end{equation}
Second we define a map $\zeta:\C^2\to\C\times\R^2$ as follows:
\[
\zeta: (z,w)\mapsto (z, r\arg(zw), r),
\]
where $r=||(z,w)||=\sqrt{|z|^2+|w|^2}$.

An explicit calculation shows that the Jacobian of this map (where defined)
is constant with value ``$-1$''.
Furthermore the gradient $\grad(r)$ is of length one and normal
on the level set $S(r)$.
Correspondingly the map $\zeta$ is volume preserving and $S(r)$
has the same volume as its image 
\begin{equation}
\label{param}
\zeta(S(r))=\{z\in\C:|z|\le r\}\times[0,2\pi r)\times\{r\},
\end{equation}
namely $2\pi^2r^3$.

Similarly the euclidean volume of $\Omega_{K,\lambda,\mu}$ agrees with the
euclidean volume of
\[
\zeta(\Omega_{K,\lambda,\mu})=\{z\in\C:|z|\le Kr^\lambda\}\times
\{\theta r:\theta\in [0,2\pi), |\sin \theta |\le r^{-\mu} \}\times\{r\}
\]
Using \eqref{sinvol}
it follows that for $r\ge 1$ the volume of $\Omega_{K,\lambda,\mu}$ is bounded
by
\[
\pi\left(Kr^\lambda\right)^2 \ \cdot \ 2r^{-\mu}\pi r = 2K^2\pi^2 r^{2\lambda+1-\mu}.
\]
In particular, 
\begin{equation}\label{eqvol}
\mathrm{vol} (\Omega_{K,\lambda,\mu})=O(r^{2\lambda+1-\mu}).
\end{equation}

%

\subsection{Arithmetic estimates.}

Besides the Landau $O$-symbols we also use the notation ``$\gtrsim$'':
If $f,g$ are functions of a real parameter $r$, then $f(r)\gtrsim g(r)$
indicates that 
\[
\liminf_{r\to +\infty}\frac{f(r)}{g(r)}\ge 1.
\]
Similarly
$f\sim g$ indicates
\[
\lim_{r\to +\infty}\frac{f(r)}{g(r)}= 1.
\]

In the sequel,
we will work with domains $\Omega\subset S(r)$
 (i.e.~for each $r>0$ some subset  $\Omega=\Omega_r\subset S(r)$
is chosen). In this context, given functions $f$, $g$ on $\C^2$ we say 
``$f(z, w)\gtrsim g(z, w)$ holds on $\Omega$'' if for every sequence
$(z_n, w_n)\in\C^2$ with $\lim ||(z_n, w_n)||=+\infty$ and
$(z_n, w_n)\in\Omega_r$ ($r=||(z_n, w_n)||$) we have
\[
\liminf_{n\to\infty}\frac{f(z_n, w_n)}{g(z_n, w_n)}\ge 1.
\]

We develop some estimates for $\phi(z, w)=|z|^2+|1+zw|^2$.
Fix $\mu>0$, $-\infty < \lambda<1$.
\begin{enumerate}
\item
For all $z, w$: $\phi\geq |z|^2$.
\item
If $(z, w)\in S(r)$ and $|z| \leq\frac{1}{2r}$, then 
\[
|w|\leq r \implies |zw|\leq\frac12 \implies |1+zw|\geq \frac12
\]
 and therefore $\phi\geq \frac14$.
\item
For $|z|\leq r^\lambda$ we have $|w|\sim r$,
i.e., for fixed $\lambda,\mu$ and any choice of $(z_r, w_r)\in S(r)$ with 
$|z_r|\leq r^\lambda$
we have $\lim_{r\to\infty}  |w_r|/r=1$.
\item
For $|z|\geq \frac{3}{2r}$ and $|z|\leq r^\lambda$ we have that
$\phi\gtrsim\frac{1}{9}|zw|^2$,
because $|w| \sim r$ and $|zw|\gtrsim \frac{3}{2}$
(equivalently, $1 \lesssim \frac{2}{3} |zw|$),
implying $|1+zw|\geq |zw| -1 \gtrsim\frac{1}{3}|zw|$.
\item
For all $z, w$,  $\phi\geq \left|\Im(1+zw)
\right|^2=\left(|zw|\sin \theta \right)^2$.
\end{enumerate}

\subsection{Putting things together.}

We are going to prove first the claim 
\[
``I(r)=O(r^{2+\varepsilon}), \quad \forall\varepsilon>0 '' 
\]
by dividing $S(r)$ into regions
$A$, $B$, $C$, $D_{-2}$, $D_{-1}$, $D_0$, $D_1$, $E$, $F$,
each of which is investigated separately.

\begin{itemize}
\item
Region $A$ consists of those points with $|z|\leq \frac1{2r}$, i.e.,
$A=\Omega_{\frac12,-1,0}$.
The volume $\mathrm{vol}(A)$ is thus of order $O(r^{-1})$ 
Due to $(ii)$ the integrand $\eta$ is bounded by $\eta|_A=O(r^2)$.
It follows that
\[
\int_A\eta \, dV \le \mathrm{vol}(A) \cdot
\sup_{(z,w) \in A}\eta(z, w)=O(r).
\]
Hence the contribution of $A$ to the integral $I_r=\int_{S(r)}\eta \, dV$
is bounded by $O(r)$.
\item
Region $B$ consists of those points with $\frac{1}{2r}\leq
|z|\leq \frac{3}{2r}$
and $|\sin \theta |< \frac{1}{r}$. Thus $B\subset \Omega_{3/2,-1,1}$.
Due to \eqref{eqvol} this implies
$\mathrm{vol} (B)=O(r^{-2})$.
For the integrand $\eta|_B$ we have the bound
$\eta|_B=O(r^4)$ (using $(i)$ and $|z|\geq\frac{1}{2r}$ ).
Hence
\[
\int_B \eta \, dV \le \mathrm{vol}(B) \cdot
\sup_{(z,w) \in B}\eta(z, w)=O(r^2) ;
\]
i.e., the contribution of $B$ to the integral $I_r$ is bounded by $O(r^2)$.
\item
Region $C$ consists of those points with
 $\frac{1}{2r}\leq |z|\leq \frac{3}{2r}$
and $|\sin \theta|>\frac{1}{r}$.
Since $|w| \sim r$, $\frac{1}{2} \lesssim |zw| \lesssim \frac32$
We take the volume-compatible parameter $\psi=r\theta$ due to
\eqref{param}.
Then $\frac{1}{r} < \left|\sin \frac{\psi}{r}\right| < \frac{\psi}{r}$,
and so $\psi>1$.
Therefore
$$
J_r:=\int_{1<\psi< 2\pi r,\, \left|\sin \frac{\psi}{r}\right|>
\frac{1}{r}} \eta\, d\psi
=\int_{1<\psi<2\pi r,\, \left|\sin \frac{\psi}{r}\right|> \frac{1}{r}}
\frac{ 2 r^2}{\left(\sin \frac{\psi}{r}\right)^2} d\psi=O(r^4).
$$
Here in fact we have that there is a constant $c>1$ such that
$$
\frac{r^4}{c}  \leq J_r \leq c r^4.
$$
Therefore it follows that
\begin{equation}
\label{lowerest}
\frac{r^2}{c'} \leq
\int_C \eta \, dV=\int_{\frac{1}{2r}\leq |z|\leq \frac{3}{2r}}\,
J_r \, \frac{i}{2} dz \wedge d\bar{z} \leq  c' r^2,
\end{equation}
where $c'$ is a positive constant.
Thus the contribution of $C$ to the integral $I_r$ is bounded by $O(r^2)$.
\item
For $\gamma\in\{-2,-1,0,1\}$ let $D_\gamma$ denote the set of those points
where $|z|\geq \frac{3}{2r}$, $|z|\leq  r^{1-\varepsilon}$
and $r^{\frac\gamma2}\leq |z|\leq  r^{\frac{\gamma+1}{2}}$.
For each $\gamma$ the integrand $\eta$ is bounded on $D_\gamma$ by
$O(r^{-\gamma})$ (due to $(iv)$),
and the volume $\mathrm{vol}(D_\gamma)$ is bounded by
$O(r^{2+\gamma})$, because $D_\gamma\subset\Omega_{1,\frac{\gamma+1}{2},0}$. 
Thus the contribution of $D_\gamma$ to the integral $I_r$
is bounded by $O(r^2)$.
\item
Let $E$ denote the region where $|z|\geq r^{1-\varepsilon}$,
$|w|\geq r^{\frac 12}$.
For the integrand we have that $\eta|_E=O(r^{2\varepsilon-1})$
(using $(iv)$). 
The  volume of $E$ is bounded by the total volume of $S(r)$,
i.e., $\mathrm{vol}(E)=O(r^3)$.
Together this shows that the contribution of $E$ to $I_r$ is bounded
by  $O(r^{2+2\varepsilon})$.
\item
Let $F$ denote the region where $|w|\leq r^{\frac 12}$. 
In analogy to $(iii)$ we have $|z| \sim r$.
With $(i)$ it follows that 
$\sup_{(z, w) \in F}\eta(z, w)=O(1)$.
On the other hand the volume of $F$
agrees with the volume of $\{(z, w)\in S(r): |z|\leq r^{\frac 12}\}$
which according to \eqref{eqvol} is bounded by $O(r^{2})$.
Together this yields that the contribution of $F$ to $I_r$
 is bounded by $O(r^2)$.
\end{itemize}

Thus we have a collection of nine regions 
($A$, $B$, $C$, $D_{-2}$, $D_{-1}$, $D_0$, $D_1$, $E$, $F$)
covering the sphere $S(r)$.
For each such region $\Omega$ we have verified
\[
\int_\Omega\eta\, dV=O(r^{2+\varepsilon}), \quad \varepsilon>0.
\]
This establishes our claim
\[
I_r=O(r^{2+\varepsilon}), \quad \varepsilon>0.
\]

Furthermore, it follows from \eqref{lowerest} that
$$
r^2=O(I_r).
$$
As a consequence, the holomorphic map $f:\C^2\to S$ induced by
$f:(z, w)\mapsto (z,1+zw)$  is of order $\rho_f=1$. \hfill {\it Q.E.D.}

\section{Problems.}

Because of the results presented above it may be interesting to
recall some problems (conjectures) from \cite{n93}, \S1.4.
An $n$-dimensional compact complex manifold
$X$ is said to be {\it unirational} if there is a surjective
meromorphic map $\phi: \pnc \to X$; in this case,
if $g: \C^n \to \pnc$ is a differentiably non-degenerate meromorphic
map with order $\rho_g<2$, then $\phi\circ g: \C^n \to X$
is differentiably non-degenerate and has order less than two.
Therefore, the rationality and the unirationality of $X$ cannot be distinguished
by the existence of a differentiably non-degenerate meromorphic map
$f: \C^n \to X$ with $\rho_f<2$.

\begin{prob}
Let $X$ be a compact k\"ahler manifold of
dimension $n$.
If there is a differentiably non-degenerate meromorphic
map $f: \C^n \to X$ with order $\rho_f<2$, is $X$ unirational?
\end{prob}

At least this is true for $\dim X \leq 2$ by Corollary \ref{1-dim}
and the Main Theorem \ref{main}.

\begin{prob}
Let $f:\C \to X$ be a non-constant entire curve
into a projective (or k\"ahler) manifold $X$.
If $\rho_f<2$, then does $X$ contain a rational curve?
\end{prob}

\bigskip

\rightline{Graduate School of Mathematical Sciences}
\rightline{The University of Tokyo}
\rightline{Komaba, Meguro,Tokyo 153-8914}
\rightline{e-mail: noguchi@ms.u-tokyo.ac.jp}
\bigskip

\rightline{Mathematisches Institut}
\rightline{Ruhr-Universit\"at Bochum}
\rightline{44780 Bochum, Germany}
\rightline{e-mail: Joerg.Winkelmann@ruhr-uni-bochum.de}

\begin{thebibliography}{99}
\setlength{\itemsep}{-3pt}
\bibitem{bpv}
W. Barth, C. Peters and A. Van de Ven,
Compact Complex Surfaces, Ergebnisse Math.\ und ihrer Grenzgebiete
{\bf 4}, Springer-Verlag, 1984.
\bibitem{ko71}
K. Kodaira, Holomorphic mappings of polydiscs into compact complex manifolds,
J. Diff.\ Geometry, {\bf 6} (1971), 33-46.
\bibitem{n93}
J. Noguchi, Some problems in value distribution and hyperbolic manifolds,
K\^oky\^uroku {\bf 819} (1993), 66-79, R.I.M.S. Kyoto University.
\bibitem{no90}
J. Noguchi and T.\ Ochiai,
Geometric Function Theory in Several Complex Variables,
Math.\ Monographs Vol.\ {\bf 80}, Amer.\ Math.\ Soc., Providence,
1990 (translated from Japanese version published from Iwanami, Tokyo 1984).
\end{thebibliography}
\end{document}